\documentclass[3p,12pt]{elsarticle}
\journal{Stochastic Processes and their Applications}

\usepackage{graphicx}
\usepackage{amssymb}
\usepackage{amsthm,amsmath}

\usepackage{hyperref}

\numberwithin{equation}{section}

\newtheorem{theorem}{Theorem}[section]

\newtheorem{corollary}[theorem]{Corollary}

\newtheorem{lemma}[theorem]{Lemma}
\theoremstyle{definition}
\newtheorem{remark}[theorem]{Remark}

\renewcommand{\Pr}{\mathbb{P}}
\newcommand{\E}{\mathbb{E}}
\renewcommand{\eta}{X}

\newcommand{\Z}{E}
\newcommand{\BZ}{\mathcal{B}(E)}
\newcommand{\st}{\mathbb{R}^d}
\newcommand{\Ro}{\overline{\mathbb R}_0^d}

\newcommand{\sgn}{\mathrm{sgn}}

\newcommand{\F}{\mathcal{F}}
\newcommand{\A}{\mathcal{A}}
\newcommand{\B}{\mathcal{B}}

\newcommand{\M}{M}
\newcommand{\ME}{\M_p((0,\infty)\times \Ro)}
\newcommand{\MEE}{\M_p(\Ro)}

\newcommand{\Do}{\mathbb{D}([0,1])}
\newcommand{\D}{\mathbb{D}([0,\infty),\mathbb{R}^d)}

\newcommand{\DT}{\mathbb{D}([0,T],\mathbb{R}^d)}
\newcommand{\ob}{Z_1}

\newcommand{\Q}{Q}
\newcommand{\Xa}{X}
\newcommand{\Pia}{\Pi}
\newcommand{\Na}{N}

\begin{document}
\begin{frontmatter}

\title{Convergence to L\'evy stable processes under some weak dependence conditions}

\author{Marta Tyran-Kami\'nska}%

\ead{mtyran@us.edu.pl}
\address{Institute of Mathematics, University of Silesia, Bankowa 14, 40-007 Katowice,
Poland}

\begin{abstract}
For a strictly stationary sequence of random vectors in $\mathbb{R}^d$ we study convergence of partial sum
processes to a L\'evy stable process in the Skorohod space with $J_1$-topology. We identify necessary and
sufficient conditions for such convergence and provide sufficient conditions when the stationary sequence is
strongly mixing.
\end{abstract}

\begin{keyword}
L\'evy stable processes\sep Poisson point processes\sep functional limit
theorem\sep Skorohod topology\sep strong mixing

\MSC 60F17\sep 60F05\sep 60G51 \sep60G55
\end{keyword}

\end{frontmatter}
\section{Introduction}

Let $\{Z_j:j\ge 1\}$ be a  strictly stationary sequence of $\st$-valued random
vectors defined on a probability space $(\Omega,\F,\Pr)$. If the $Z_j$ are
i.i.d.~then according to Rva{\v{c}}eva~\cite{rvaceva62}, there exist sequences of
constants $b_n> 0$  and $c_n$  such that
\begin{equation}\label{d:iid}
\frac{1}{b_n}\sum_{j=1}^n Z_j-c_n\xrightarrow{d}\zeta_\alpha \quad \text{in }\st
\end{equation}
for some non-degenerate  $\alpha$-stable random vector $\zeta_\alpha$ with $\alpha\in(0,2)$ if and only if $Z_1$
is \emph{regularly varying with index} $\alpha\in(0,2)$: there exists a probability measure $\sigma$ on
$\B(\mathbb{S}^{d-1})$, the Borel $\sigma$-algebra of the unit sphere $\mathbb{S}^{d-1}=\{x\in\st\colon
|x|=1\}$, such that as $x\to\infty$
\begin{equation}\label{e:rv2}
\dfrac{\Pr(|Z_1|>rx, Z_1/|Z_1|\in A )}{\Pr(|Z_1|>x)}\xrightarrow{}
r^{-\alpha}\sigma(A)
\end{equation}
for all  $r>0$ and $A\in\B(\mathbb{S}^{d-1})$ such that $\sigma(\partial A)=0$ (the notation `$\xrightarrow{d}$
in $\mathbb{X}$' refers to weak convergence of distributions of given random elements with values in the space
$\mathbb{X}$ and $|\cdot|$ denotes the Euclidean norm). The sequences $c_n$ and $b_n$ can be chosen as
\begin{equation}\label{d:bncn}
 c_n=\frac{n}{b_n}\E(Z_1I(|Z_1|\le
b_n))\quad \text{and}\quad n\Pr(|Z_1| > b_n) \to 1.
\end{equation}
A result of  Ibragimov~\cite{ibragimov62}, and its extension  to random vectors by Philipp~\cite{phillipp80},
tells us that if for a strongly mixing sequence the normalized partial sums in~\eqref{d:iid} converge to a
non-degenerate random vector $\zeta$, then necessarily $\zeta=\zeta_\alpha$ for some $\alpha\in (0,2]$; the case
of $\alpha=2$ refers to a Gaussian distributed random vector.

For the functional generalization of \eqref{d:iid} define the partial sum
processes
\begin{equation}\label{eq:psum}
X_n(t)=\frac{1}{b_n}\sum_{1\le j\le nt} Z_j -tc_n,\quad t\ge 0.
\end{equation}
For each $\omega$, $X_n(\cdot)$ is an element of the Skorohod space $\D$ of all $\st$-valued functions on
$[0,\infty)$ that have finite left-hand limits and are continuous from the right. In this paper we study weak
convergence of distributions of the partial sum processes in $\D$ with the Skorohod $J_1$ topology (see
Section~\ref{s:skor}). In the i.i.d.~case if
\begin{equation}\label{eq:convst}
X_n\xrightarrow{d} X \quad \text{in } \D, \end{equation} then necessarily $X$ is a
L\'evy $\alpha$-stable process with $\alpha\in(0,2]$, whose increments are
stationary, independent, and $X(1)$ has the same distribution as~$\zeta_\alpha$.

In the case of $\alpha=2$ in~\eqref{d:iid}, a substantial amount of work has been devoted to extend the central
limit theorem and the Donsker's invariance principle to weakly dependent random variables; see the recent review
by Merlev{\`e}de, Peligrad, and Utev~\cite{merlevedeetal06} for sequences of random variables with finite
variances and Bradley~\cite{bradley88}, Shao~\cite{shao93} for random variables with infinite variances. In the
case $\alpha\in(0,2)$ it was shown by Avram and Taqqu~\cite{avramtaqqu92} that for some $m$-dependent random
variables weak--$J_1$ convergence cannot hold despite the fact that~\eqref{d:iid} holds. However, for
$\varphi$-mixing sequences there is a characterization of convergence in~\eqref{eq:convst} by
Samur~\cite{samur87} in terms of convergence in \eqref{d:iid} and some additional conditions.
Both~\cite{avramtaqqu92} and~\cite{samur87}  use the Skorohod approach~\cite{skorohod57} in $\Do$ via tightness
plus convergence of finite dimensional distributions.

Since the case of $\alpha=2$ is to some extent well understood we shall focus on
the less studied case of $\alpha<2$. Our method of proof of~\eqref{eq:convst} is
based on point process techniques used by Durrett and
Resnick~\cite{durrettresnick78} for convergence of dependent random variables. For
a comprehensive account on this subject in the independent case, we refer the
reader to the expository article by Resnick~\cite{resnick86} and to his recent
monograph~\cite{resnick07}. We recall relevant notation and background in
Section~\ref{s:pr}. In Section~\ref{s:nsc} we study, in a somewhat more general
setting, the problem of convergence $X_n\xrightarrow{d}X$ in $\D$ with $X$ being a
L\'evy process without Gaussian component. Theorem~\ref{th:main} gives necessary
and sufficient conditions for such convergence. Roughly speaking these are a
convergence of point processes $N_n$, consisting of the jump points of $X_n$, to
the corresponding point process $N$ of jumps of $X$, which necessarily is a
Poisson point process, and a condition which allows one to neglect accumulation of
small jumps. In the case of a strictly stationary strongly mixing sequence
$\{Z_j\colon j\ge 1\}$  of random variables, the class of possible limiting
processes for $N_n$ were discussed first by Mori~\cite{mori77}. Any element of
this class must  be infinitely divisible and invariant under certain
transformations. Then Hsing~\cite{hsing87} derived cluster representations of the
limiting processes under somewhat weaker distributional mixing conditions. Thus we
need to impose extra condition to obtain a Poisson process in the limit. We now
describe an application of Theorem~\ref{th:main} to a strictly stationary sequence
$\{Z_j\colon j\ge 1 \}$ under strong mixing conditions.

There exist several coefficients 'measuring' the dependence between two
$\sigma$-algebras $\A$ and $\B\subset \F$, the most usual ones being
\[
\alpha(\A,\B)=\sup\{|\Pr(A\cap B)-\Pr(A)\Pr(B)|: A\in\A, B\in\B\},
\]
\[
\varphi(\A,\B)=\sup\{|\Pr(B|A)-\Pr(B)|: A\in\A, \Pr(A)>0, B\in\B\},
\]
and the maximal coefficient of correlations
\[
\rho(\A,\B)=\sup\{|\mathrm{Corr}(f,g)|: f\in L^2(\A),g\in L^2(\B)\Bigr\};
\]
see the review paper by Bradley~\cite{bradley05} as a general reference for mixing
conditions.

Given the sequence $\{Z_j\colon j\ge 1\}$, we define $\F_{m}^n=\sigma\{Z_j:m\le
j\le n\}$ and, for every $ n\ge 1$,
\[
\phi_0(n)=\sup_{k\ge 1}\alpha(\F_{1}^k, \F_{n+k}^\infty),
\]
\[
\phi_1(n)=\sup_{k\ge 1}\varphi(\F_{1}^k, \F_{n+k}^\infty),
\]
\[
\rho(n)=\sup_{k\ge 1}\rho(\F_{1}^k, \F_{n+k}^\infty).
\]
The sequence $\{Z_j\colon j\ge 1\}$  is said to be \emph{mixing with rate function
$\phi_s$} if $\phi_s(n)\to 0$  as $n\to\infty$; the case of $s=0$ ($s=1$) refers
to \emph{strongly (uniformly} or \emph{$\varphi$-) mixing} sequence.

Our main result for strongly mixing sequences is the following functional limit
theorem.
\begin{theorem}\label{t:rho}
Let a strictly stationary sequence $\{Z_j:j\ge 1\}$ be mixing with rate function~$\phi_0$. Assume that $Z_1$ is
regularly varying with index $\alpha\in(0,2)$ and that one of the following conditions is satisfied:
\begin{enumerate}
\item\label{l:rho1} $\alpha\in (0,1)$;
\item\label{l:rho2} $\alpha\in [1,2)$ and for every $\delta>0$
\begin{equation}\label{eq:esmall1d1}
\lim_{\varepsilon\to 0}\limsup_{n\to\infty}\Pr\bigl(\max_{1\le k\le n}|\sum_{j=1}^{k}(Z_{j}I(|Z_{j}|\le
\varepsilon b_n)-\E(Z_{1}I(|Z_{1}|\le \varepsilon b_n)))|\ge \delta b_n\bigr)= 0.
\end{equation}
\end{enumerate}
Then $X_n\xrightarrow{d}X$ in $\D$, where $X_n$ is as in
\eqref{d:bncn}--\eqref{eq:psum}, and  $X$ is a L\'evy $\alpha$-stable process if
and only if  the following local dependence condition holds:
\begin{description}
\item[\textbf{LD}$(\phi_0)$] For any $\varepsilon>0$ there exist sequences of integers
$r_n=r_n(\varepsilon),l_n=l_n(\varepsilon)\to\infty$ such that
\begin{equation}\label{eq:asneg0}
r_n=o(n),\quad l_n=o(r_n),\quad n\phi_0(l_n)=o(r_n), \quad \text{as }n\to \infty,
\end{equation}
 and
\begin{equation}\label{eq:asneg}
\lim_{n\to\infty}\Pr(\max_{2\le j\le r_n}|Z_j|> \varepsilon b_n
\bigl||Z_1|>\varepsilon b_n)=0.
\end{equation}
\end{description}
\end{theorem}

 Note
that if $\phi_0(n)\to 0$ then one can always find sequences $r_n$, $l_n$ satisfying \eqref{eq:asneg0}. Thus, if
for any $\varepsilon>0$
\begin{equation*}
\lim_{k\to\infty}\limsup_{n\to\infty}\Pr(\max_{2\le j\le [n/k]}|Z_j|> \varepsilon b_n \bigl| |Z_1|>\varepsilon
b_n)=0,
\end{equation*}
then condition \textbf{LD}$(\phi_0)$ holds. Therefore, \textbf{LD}$(\phi_0)$ is implied by \eqref{eq:D'} which
is the local dependence condition $D'$ of Davis~\cite{davis83}, since $n\Pr(|Z_1|>\varepsilon b_n)\to
\varepsilon^{-\alpha}$ by~\eqref{e:rv2} and~\eqref{d:bncn}.

\begin{corollary}\label{c:rho}
Let a strictly stationary sequence $\{Z_j:j\ge 1\}$ be mixing with rate function~$\phi_0$. Assume that $Z_1$ is
regularly varying with index $\alpha\in(0,2)$. If for any $\varepsilon>0$
\begin{equation}\label{eq:D'}
\lim_{k\to\infty}\limsup_{n\to\infty}n\sum_{j=2}^{[n/k]}\Pr(|Z_j|> \varepsilon b_n ,|Z_1|>\varepsilon b_n)=0,
\end{equation}
then condition \textbf{LD}$(\phi_0)$ holds.
\end{corollary}

For uniformly mixing sequences we have the following result.
\begin{corollary}\label{c:kobus}
Let a strictly stationary sequence $\{Z_j:j\ge 1\}$ be mixing with rate
function~$\phi_1$. Assume that $Z_1$ is regularly varying. Then condition
\textbf{LD}$(\phi_0)$ is equivalent to
\begin{description}
\item[\textbf{LD}$(\phi_1)$] For any $\varepsilon>0$ and $j\ge2$
\begin{equation}\label{d:watson} \lim_{n\to\infty}\Pr(|Z_{j}|>\varepsilon
b_n\bigl| |Z_1|>\varepsilon b_n)=0.
\end{equation}
\end{description}
\end{corollary}

With the notation as in Theorem~\ref{t:rho} we have the following
characterization for $m$-dependent sequences.

\begin{corollary}\label{c:mdep}
Assume that $\{Z_j:j\ge 1\}$ is $m$-dependent. Then $X_n\xrightarrow{d} \Xa$ in
$\D$ if and only if $Z_1$ is regularly varying and for any $\varepsilon>0$ and $j=
2,\ldots,m$ condition \eqref{d:watson} is satisfied.
\end{corollary}

The proofs of these results are presented in Section~\ref{s:sms}. The reader is referred to
Theorem~\ref{t:stable} and Remark~\ref{r:asneg} for necessity of conditions~\eqref{eq:esmall1d1}
and~\eqref{eq:asneg} without the strong mixing assumption and to~Remark~\ref{r:exi} for a relation
between~\eqref{eq:asneg}~and the extremal index of the sequence $\{|Z_j|\colon j\ge 1\}$. We show in
Lemma~\ref{lemat} that condition \eqref{l:rho2} of Theorem~\ref{t:rho} holds if $\sum_{j}\rho(2^j)<\infty$.
Thus, Theorem~\ref{t:rho} complements the results of Bradley~\cite{bradley88} and Shao~\cite{shao93}.
Theorem~\ref{t:rho} together with Corollary~\ref{c:kobus} establishes Corollary 5.9 of Kobus~\cite{kobus95},
which was proved using the results of Samur~\cite{samur87}.

The methods and results of this paper were used in~\cite{tyran09b} to prove functional limit theorems for
particular examples of stationary sequences arising from dynamical systems such as continued fractions,
Gibbs-Markov maps, and piecewise expanding
maps~\cite{aaronson86,aaronsonnakada05,iosifescukraaikamp02,melbournenicol09,samur89,zweimuller}. In that
setting condition~\eqref{eq:asneg} has a nice interpretation in terms of hitting times and it can be also used
without the strong mixing assumption, see~\cite[Sections 3 and 4]{tyran09b}.

We should also point out that proving the weak convergence of the partial sums of a strictly stationary sequence
to an infinite variance $\alpha$-stable random vector in~\eqref{d:iid} might require less restrictive
assumptions as opposed to the weak convergence in~\eqref{eq:convst}. The recent paper~\cite{bartkiewicz09}
contains a detailed study of sufficient conditions for the convergence in~\eqref{d:iid} and a comparison of
various conditions used to prove stable limits under the assumption that the stationary sequence is jointly
regularly varying, which means that all finite dimensional distributions are regularly varying with index
$\alpha\in(0,2)$. Here we comment on the approach through point processes. A number of authors studied the point
processes $N_n'$ consisting of the points $b_n^{-1} Z_j$, $j=1,\ldots,n$, in order to obtain convergence to
$\alpha$-stable random vectors. The one dimensional case with $N_n'$ converging to a Poisson process was studied
in~\cite{davis83}. A systematic application of point process techniques for obtaining limit theorems for arrays
of dependent random vectors has been developed in~\cite{jakubowskikobus89,kobus95}. Their results for
 strongly mixing stationary sequences are obtained under stronger
assumptions than our condition \textbf{LD}$(\phi_0)$.  For a jointly  regularly varying  stationary sequence of
dependent random variables sufficient conditions for convergence in \eqref{d:iid} can be found in
\cite{davishsing95} and their multivariate extensions in \cite{davismikosch98}; here the limiting process for
$N_n'$ might not be a Poisson process, so that their examples provide a large class of processes for which the
functional limit theorem does not hold in the $J_1$-topology. It would be interesting to obtain corresponding
results in one of the weaker Skorohod's topologies $M_1$ or $M_2$ as defined in \cite{skorohod56};
see~\cite{avramtaqqu92} for a result in the $M_1$ topology.  After the submission of this paper, we became aware
of~\cite{basrak10}, where sufficient conditions for $X_n\xrightarrow{d} X$ in $\mathbb{D}([0,1],\mathbb{R})$
with the $M_1$ topology were obtained by building upon the approach and assumptions in~\cite{davishsing95}; one
 condition is the same as our condition~\eqref{l:rho2} in Theorem~\ref{t:rho}.

\section{Preliminaries}\label{s:pr}

In this section we collect some basic tools and notions to be used throughout this
paper.
\subsection{Point  processes}
We begin by introducing some background on point processes. We follow the point
process theory as presented in Kallenberg~\cite{kallenberg75} and
Resnick~\cite{resnick87}. Let $\Z$ be a locally compact Hausdorff topological
space with a countable basis for its topology. For our purposes, $\Z$ is a subset
of either $\Ro:=\overline{\mathbb{R}}^d\setminus\{0\}$ or $[0,\infty)\times \Ro$,
where $\overline{\mathbb{R}}=\mathbb{R}\cup\{-\infty,\infty\}$. The topology on
$\Ro$ is chosen so that the Borel $\sigma$-algebras $\B(\Ro)$ and $\B(\st)$
coincide on $\st\setminus\{0\}$. Moreover, $B\subset\Ro$ is relatively compact (or
bounded) if and only if $B\cap \st$ is bounded away from zero in $\st$, i.e.,
$0\notin \overline{B\cap \st}$.

Let $\M(\Z)$ be the set of all
 Radon measures on $\BZ$, i.e., nonnegative Borel measures which are finite
on relatively compact subsets of $\Z$. The space $\M(\Z)$ is a Polish space when
considered with the topology of vague convergence. Recall that $m_n$
\emph{converges vaguely} to $m$
\[
m_n\xrightarrow{v} m\quad \text{iff}\quad  m_n(f)\to m(f) \quad \text{for all }
f\in C^+_K(\Z),
\]
where $m(f)=\int_{\Z} f(x) m(dx)$ and $C^+_K(\Z)$ is the space of nonnegative continuous functions on $\Z$ with
compact support. We have $m_n\xrightarrow{v} m$ if and only if $m_n(B)\to m(B)$ for all relatively compact $B$
for which $m(\partial B)=0$.

The set $\M_p(\Z)$ of point measures on $\Z$ is a closed subspace of $\M(\Z)$
consisting of all integer-valued measures in $\M(\Z)$. Denote by $\epsilon_x$ the
unit measure concentrated at $x\in \Z$. Any point measure $m\in \M_p(\Z)$ is of
the form $m=\sum_i\epsilon_{x_i}$, where $\{x_i\}$ is at most a countable
collection of points of $\Z$. The measure $m$ is called \emph{simple} if the
points $\{x_i\}$ are all distinct.

A point process $N$ on $\Z$ is an $\M_p(\Z)$-valued random variable, defined on a
probability space $(\Omega ,\mathcal{F}, \Pr)$.  The measure $\Q$ defined by
$\Q(A)=\E (N(A))$, $A\in\B(\Z)$, is called a \emph{mean measure of $N$}. The
process $N$ is called \emph{simple} if almost all its realizations are simple. A
point process $N$ is called a \emph{Poisson process} with mean measure $\Q\in
\M(\Z)$ if $N(A_1),\ldots,N(A_l)$ are independent random variables for any
disjoint sets $A_1,\ldots,A_l\in \B(\Z)$ and $N(A)$ is a Poisson random variable
with mean $\Q(A)$ for $A\in \B(\Z)$ with $\Q(A)<\infty$.  The Poisson process is
simple if its mean measure is non-atomic. The Laplace functional of Poisson
process $N$ is of the form
\[
\E  [e^{-N(f)} ]= \exp\left\{-\int_{\Z} (1 - e^{-f(x)})\Q(dx)\right\}
\]
for nonnegative measurable $f$. Given a sequence of point processes $N_n$ we have
$N_n\xrightarrow{d} N$ in $\M_p(\Z)$, by~\cite[Theorem 4.2]{kallenberg75}, if and
only if $\E [e^{-N_n(f)} ]\to \E [e^{-N(f)}]$ for all $ f\in C^+_K(\Z)$.

\subsection{Infinitely divisible and stable random vectors}

An infinitely divisible $\st$-valued random vector $\zeta$ is uniquely determined
through the L\'evy-Khintchine formula, which states that its characteristic
function is of the form
\begin{equation*}\label{e:infdiv1}
\E e^{i\langle u, \zeta\rangle}= \exp\bigl(i\langle a,u\rangle -\frac12\langle
\Sigma u,u\rangle + \int_{\st}\bigl( e^{i\langle u,x\rangle}-1 -i\langle
u,x\rangle I(|x|\le 1)\bigr)\Pi(dx)\bigr),
\end{equation*}
where $a$ is a $d$-dimensional vector, $\Sigma$ is a  symmetric nonnegative
definite $d \times d$ matrix, and $\Pi$ is a \emph{L\'evy measure}, i.e., a
$\sigma$-finite Borel measure on $\st$ such that $\Pi(\{0\})=0$ and
\begin{equation*}
\int_{\st} (1\wedge |x|^2) \Pi(dx) <\infty.
\end{equation*}
Here $\langle x,u\rangle$  denotes the usual inner product in $\st$. We have $\Pi(B)<\infty$ for any set
$B\in\B(\st)$ bounded away from $0$. We can extend $\Pi$ on $\B(\Ro)$ in such a way that
$\Pi(\overline{\mathbb{R}}^d\setminus \st)=0$. The generating triplet $(\Sigma,\Pi,a)$ uniquely determines a
given infinitely divisible random vector.

A particular class of infinitely divisible random vectors without Gaussian component, i.e., with $\Sigma=0$, are
stable random vectors. Recall that a vector $\zeta$ with values in $\st$ is said to be $\alpha$-stable for some
$\alpha\in (0, 2)$ if its characteristic function is given by
\[
\E e^{i\langle u, \zeta\rangle}=
\left\{
\begin{array}{ll}
\exp(-\int_{\mathbb{S}^{d-1}}|\langle
u,s\rangle|^\alpha(1-i\tan\frac{\pi\alpha}{2}\sgn\langle
u,s\rangle)\lambda_1(ds)+i\langle \tau,u\rangle)
 , &   \alpha\neq 1,\\
\exp(-\int_{\mathbb{S}^{d-1}}|\langle u,s\rangle|(1+i\frac{2}{\pi}\sgn\langle
u,s\rangle\log|\langle u,s\rangle|)\lambda_1(ds)+i\langle \tau,u\rangle) , &
\alpha= 1,
\end{array}
\right.
\]
where $\tau\in\st$ and  $\lambda_1$ is a finite nonzero measure on
$\B(\mathbb{S}^{d-1})$; see~\cite[Theorem 14.3]{sato99}. Its L\'evy measure is of
the form
\begin{equation}\label{d:Pia}
\Pi_\alpha(B)=\int_{S^{d-1}}\int_{0}^\infty
1_B(rs)r^{-\alpha-1}dr\lambda(ds),\quad B\in\B(\mathbb{R}^d),
\end{equation} where
$\lambda$ is a finite nonzero measure on  $\B(\mathbb{S}^{d-1})$ and it is a constant multiple of $\lambda_1$.

\subsection{Regularly varying vectors}
If $\zeta$ is $\alpha$-stable for some $\alpha\in (0, 2)$, it is regularly varying
with index $\alpha$ as defined in Introduction. The concept of vague convergence
in $\M(\mathbb{S}^{d-1})$ allows us to rewrite condition~\eqref{e:rv2} in the form
\begin{equation}\label{e:rv3}
\dfrac{\Pr(|\zeta|>tx, \zeta/|\zeta|\in\cdot )}{\Pr(|\zeta|>x)}\xrightarrow{v}
t^{-\alpha}\sigma(\cdot)\quad \text{as}\quad x\to\infty.
\end{equation}
The measure  $\sigma$ is called the \emph{spectral measure}. See \cite[Chapter
5]{resnick87} for background on multivariate regular variation.

Alternatively, \eqref{e:rv3} is equivalent to: there exists $\Q\in M(\Ro)$ such that
$\Q(\overline{\mathbb{R}}^d\setminus\st)=0$ and
\begin{equation*}
\dfrac{\Pr(x^{-1} \zeta \in\cdot )}{\Pr(|\zeta|>x)}\xrightarrow{v} \Q\quad
\text{in } \M(\Ro), \quad \text{as}\quad x\to\infty,
\end{equation*}
and to the sequential definition of regular variation: there exists $b_n\to\infty$
such that
\begin{equation}\label{d:rv}
n\Pr(b_n^{-1}\zeta \in \cdot)\xrightarrow{v} \Q,\quad n\to\infty.
\end{equation}
One can always choose $b_n$ such that $n\Pr(|\zeta|>b_n)\sim 1$. The measure $\Q$ necessarily has the property
$\Q(r B)=r^{-\alpha}\Q(B)$, $r>0$, for some $\alpha>0$ and all $b\in\B(\Ro)$, which clarifies the relation with
the index $\alpha$.

If \eqref{d:rv} holds with $\Q=\Pi_\alpha$ where $\Pi_\alpha$ is as in
\eqref{d:Pia}, then for the finite measure $\lambda$ in \eqref{d:Pia} we have
$\lambda=\alpha \sigma$. Note that the property of regular variation does not
dependent on a given norm $|\cdot|$ in $\st$, however, the spectral measure and
the limit measure $Q$ are different for distinct norms.

Regular variation of $\zeta$ implies that the function $x\mapsto \Pr(|\zeta|>x)$
is regularly varying:
\begin{equation}\label{e:xitail}
  \Pr(|\ob|>x)=x^{-\alpha}L(x),
\end{equation}
where $L$ is a \emph{slowly varying function}, i.e., $L(rx)/L(x)\to 1$ as
$x\to\infty$ for every $r>0$.  In the special case of $d=1$  a random variable
$\zeta$ is regularly varying with index $\alpha\in(0,2)$ if and only if
\eqref{e:xitail} holds for a slowly varying $L$ and the tails are balanced: there
exist $p,q\ge 0$ with $p+q = 1$ such that
\begin{equation}\label{e:xitail1}
\lim_{x\to\infty}\dfrac{\Pr(\ob>x)}{\Pr(|\ob|>x)}=p\quad \mbox{and}\quad
\lim_{x\to\infty}\dfrac{\Pr(\ob<-x)}{\Pr(|\ob|>x)}=q.
\end{equation}
From \eqref{e:xitail1} and \eqref{e:xitail} it follows that
\begin{equation*}
\lim_{n\to\infty}\dfrac{nL(b_n)}{b_n^\alpha}=1
\end{equation*}
and for $r>0$ we have
\[
\lim_{n\to\infty}n\Pr(Z_1>b_n
r)=r^{-\alpha}p\quad\mbox{and}\quad\lim_{n\to\infty}n\Pr(Z_1<-b_n r)=r^{-\alpha}q.
\]
Hence
\begin{equation*}
    n\Pr(b_n^{-1}Z_1\in \cdot )\xrightarrow{v} \Pi_\alpha(\cdot),
\end{equation*}
where $\Pi_\alpha$ is an absolutely continuous measure on $\mathbb{R}$ with density
\begin{equation*}
  \Pi_\alpha(dx)=\left(p\alpha 1_{(0,\infty)}(x)+q\alpha
1_{(-\infty,0)}(x)\right)|x|^{-\alpha-1}dx.
\end{equation*}

\subsection{Skorohod $J_1$ topology}\label{s:skor}
Let $\DT$ for $T>0$ be the space of all $\st$-valued functions $\psi$ on $[0,T]$ that are right continuous on
$[0,T)$ and have finite left-hand limits $\psi(t-)$ for all $t\in(0,T]$. Consider the set $\Lambda_T$ of
strictly increasing, continuous mappings $\lambda$ of $[0,T]$ onto itself such that $\lambda(0)=0$ and
$\lambda(T)=T$. The Skorohod $J_1$ metric on $\DT$ is defined as
\[
d_{T}(\psi_1,\psi_2)=\inf_{\lambda\in \Lambda_T}(\sup_{0\le s\le T}|\psi_1(\lambda(s))-\psi_2(s)|\vee \sup_{0\le
s\le T}|\lambda(s)-s|)
\]
for $\psi_1,\psi_2\in\DT$, where $a\vee b=\max\{a,b\}$. Then
\[
d_\infty(\psi_1,\psi_2)=\int_{0}^\infty e^{-t}(d_t(\psi_1,\psi_2)\wedge 1)dt,\quad \psi_1,\psi_2\in\D,
\]
defines a metric in $\D$ which induces the Skorohod $J_1$-topology. The metric spaces $(\DT,d_T)$ and
$(\D,d_\infty)$ are separable. For more details see \cite{whitt80} and \cite[Section 6]{jacodshiryaev03}.

We have the following characterization of weak convergence in $\D$. If $X_n, X$ are stochastic processes with
sample paths in $\D$ then $X_n\xrightarrow{d} X$ in $\D$ if and only if $X_n\xrightarrow{d} X$ in $\DT$ for all
$T\in T_X=\{t>0:\Pr(X(t)\neq X(t-))=0\}$.

\section{Necessary and sufficient conditions for weak convergence in the Skorohod topology}
\label{s:nsc}

In this section we study convergence in distributions of partial sum processes to L\'evy processes without
Gaussian component in the Skorohod space $\D$ with $J_1$-topology.

 Let $\eta=\{\eta(t):t\ge 0\}$  be an $\st$-valued L\'evy process, i.e., $\eta(0)=0$ a.s,
$\eta$ has stationary independent increments and sample paths in $\D$. We assume that $X$ is such that
\begin{equation}\label{eq:le1}
\E (e^{i\langle u, \eta(1)\rangle})=\exp \bigl(\int \bigl(e^{i\langle u,
x\rangle}-1-i\langle u, x\rangle I(|x|\le 1)\bigr)\Pi(dx)\bigr),
\end{equation}
where  $\Pi$ is a L\'evy measure with $\Pi(\{x:|x|=1\})=0$.  The jump process
$\Delta \eta(t):=\eta(t)-\eta(t-)$, $t>0$, determines a Poisson point process $N$
on $[0,\infty)\times\Ro$, which can be represented as
\begin{equation}\label{d:pp}
N=\sum_{\{t:\Delta\eta(t)\neq 0\}}\epsilon_{(t,\Delta\eta(t))},
\end{equation}
and  $\mathrm{Leb}\times \Pi$ is the mean measure of $N$; see \cite{sato99} for
details.

We also consider a family  $\{\eta_{n,j}\colon j, n\ge 1\}$ of $\st$-valued random
vectors such that $\Pr(0<|\eta_{n,j}|<\infty)=1$. Define the partial sum process
\begin{equation*}\label{d:xn}
X_n(t)=\sum_{1\le j\le nt}\eta_{n,j}-tc_n,\quad t\ge 0,n\ge 1,
\end{equation*}
where  $\{c_n\colon n\ge 1\}$ is a sequence of vectors in $\st$. It should be emphasized that we are not
assuming any dependence structure for the random vectors $\{\eta_{n,j}\colon j, n\ge 1\}$ in this section.

The following result extends \cite[Theorem 2.10.1]{leadbetterrootzen88},
\cite[Theorem 4.1]{durrettresnick78}, \cite[Proposition 3.4]{resnick86}. It
describes the connection between convergence in the Skorohod space with
$J_1$-topology and convergence of the corresponding point processes of jumps in
$\ME$.

\begin{theorem}\label{th:main} Let $X$ be a
L\'evy process satisfying \eqref{eq:le1} and let $N$ be the corresponding Poisson point process as
in~\eqref{d:pp}. Then
\begin{equation}\label{d:con}
X_n\xrightarrow{d} X \quad \text{in }\D
\end{equation}
if and only if
\begin{equation}\label{eq:pce}
N_n:=\sum_{j\ge 1}\epsilon_{(\frac{j}{n},\eta_{n,j})}\xrightarrow{d} N \quad
\text{in }\ME
\end{equation}
and for every $\delta>0$ and $T>0$
\begin{equation}\label{ce:assmall}
\lim_{\varepsilon\to 0}\limsup_{n\to\infty}\Pr\bigl(\sup_{0\le t\le T}|\sum_{j\le
nt}\eta_{n,j}I(|\eta_{n,j}|\le \varepsilon)-t(c_n-\int_{\{x:\varepsilon <|x|\le
1\} }x\Pi(dx))|\ge \delta\bigr)= 0,
\end{equation}
where the limit is taken over all $\varepsilon$ with
$\Pi(\{x:|x|=\varepsilon\})=0$.
\end{theorem}

We shall prove this result in Section~\ref{s:proof} using the continuous mapping theorem and properties of
L\'evy processes. In the next section we will use this result to prove limit theorems for stationary sequences.
For completeness we also provide the following sufficient conditions for convergence of marginal distributions.
\begin{theorem}\label{p:marg}
Let $N'$ be a Poisson point process in $\MEE$ with mean measure $\Pi$. If
\begin{equation}\label{eq:pce1d}
N_n':=\sum_{j=1}^n\epsilon_{\eta_{n,j}}\xrightarrow{d} N' \quad \text{in }\MEE
\end{equation}
and for any $\delta>0$
\begin{equation}\label{ce:assmall1d}
\lim_{\varepsilon\to 0}\limsup_{n\to\infty}\Pr\bigl(|\sum_{j=1}^n
\eta_{n,j}I(|\eta_{n,j}|\le \varepsilon)-c_n+\int_{\{x:\varepsilon <|x|\le 1\}
}x\Pi(dx)|\ge \delta\bigr)= 0,
\end{equation}
where the limit is taken over all $\varepsilon$ with
$\Pi(\{x:|x|=\varepsilon\})=0$, then
\begin{equation}\label{eq:cmarg}
X_n(1)\xrightarrow{d} X(1)\quad \text{in }\st.
\end{equation}
\end{theorem}

\begin{remark}
Note that Poisson convergence in \eqref{eq:pce1d} is not necessary for~\eqref{eq:cmarg}. There are many examples
of dependent random variables for which the latter holds, but in \eqref{eq:pce1d} we have convergence to a
non-Poisson point process \cite{davishsing95}.

Theorem~\ref{th:main} can be used to disprove weak convergence in Skorohod $J_1$-topology. Since
condition~\eqref{eq:pce} implies~\eqref{eq:pce1d}, the convergence in \eqref{d:con} is impossible in the
$J_1$-topology for all examples where $N'$ was shown to be non-Poisson. These include moving average processes
\cite{davisresnick85} for which lack of convergence in the $J_1$-topology was shown in \cite{avramtaqqu92} using
the finite dimensional plus tightness technique.
\end{remark}

We conclude this section with a discussion of some of our conditions. A situation
when condition \eqref{ce:assmall} is not needed at all is described in the
following.

\begin{corollary}
Let the L\'evy measure $\Pi$ be such that $\int (1\wedge|x|)\Pi(dx)<\infty$ and
let $c_n=0$ for all $n\ge 1$. Suppose that for any $T>0$ we have
\begin{equation}\label{e:fm}
\lim_{\varepsilon\to 0}\limsup_{n\to\infty}\sum_{j\le
nT}\E(|\eta_{n,j}|I(|\eta_{n,j}|\le \varepsilon))=0.
\end{equation}
Then $N_n\xrightarrow{d} N$ in $\ME$ if and only if $X_n\xrightarrow{d} \tilde{X}$
in $\D$, where
\[
\tilde{X}(t)=X(t)+t\int xI(|x|\le 1)\Pi(dx),\quad t\ge 0.
\]
\end{corollary}

This result is a consequence of Theorem~\ref{th:main} and the following maximal inequality from~\cite[Theorem
1]{kounias}, the proof of which extends directly to random vectors. Note that if $\Pi=\Pi_\alpha$, where
$\Pi_\alpha$ is as in \eqref{d:Pia} with $\alpha\in (0,2)$, then  $\int (1\wedge|x|)\Pi(dx)<\infty$ precisely
when $\alpha<1$.
\begin{lemma}\label{l:max}
If $\zeta_i$ are $\st$-valued random variables with $E|\zeta_i|<\infty$,
$i=1,\ldots,k$, then for any $\delta>0$ we have
\[
\Pr(\max_{1\le j\le k}|\zeta_1+\ldots+\zeta_j|\ge \delta)\le
\frac{1}{\delta}\sum_{i=1}^k \E|\zeta_i|.
\]
\end{lemma}

\begin{remark}\label{l:iddis}
Condition \eqref{ce:assmall} can be equivalently replaced by
\begin{equation}\label{eq:esmall1d}
\lim_{\varepsilon\to 0}\limsup_{n\to\infty}\Pr\bigl(\sup_{0\le t\le T}|\sum_{j\le
nt}(\eta_{n,j}I(|\eta_{n,j}|\le \varepsilon)-\E(\eta_{n,j}I(|\eta_{n,j}|\le \varepsilon)))|\ge \delta\bigr)= 0
\end{equation}
if for example, for any $\varepsilon>0$,
\begin{equation*}
\lim_{n\to\infty}\sup_{0\le
t\le T}|\sum_{j\le nt}\E(X_{n,j} I(|X_{n,j}|\le \varepsilon
))-t(c_n-\int_{\{x\colon\varepsilon<|x|\le 1\}}x\Pi(dx))|=0.
\end{equation*}
\end{remark}

In order to check condition~\eqref{eq:pce} we may use Kallenberg's characterization of convergence to simple
point processes. In particular, the following lemma  is a direct consequence of ~\cite[Theorems 4.7 and
4.8]{kallenberg75} with later improvements from \cite{kallenberg96}.

\begin{lemma}\label{l:pc}
Let $N$ be a Poisson point process with mean measure $\mathrm{Leb}\times \Pi$ where the L\'evy measure $\Pi$ is
non-atomic and let
 $\mathcal{U}$ be the class of all finite unions of rectangles in
$\Ro$ bounded away from $0$ and with boundary of zero $\Pi$ measure.
\begin{enumerate}
\item\label{l:pcl1} We have
\[
N_n\xrightarrow{d} N \quad \text{in } \ME
\]
if and only if
\begin{equation}\label{eq:1dpc}
N_n((s,t]\times \cdot) \xrightarrow{d} N((s,t]\times \cdot) \quad \text{in } \MEE
\end{equation}
for all $t>s\ge 0$ and
\begin{equation*}
\lim_{n\to\infty}\E(e^{-N_n(B)})=\E(e^{-N(B)})
\end{equation*}
for any set $B$ of the form $\bigcup_{j=1}^k (s_j,t_j]\times U_j$, where $0\le
s_1<t_1<\ldots <s_k<t_k$, and $U_j\in \mathcal{U}$, $j=1,\ldots,k$, $k\ge 2$.
\item\label{l:pcl2} Let  $t>s\ge
0$. If $\E N_n((s,t]\times \cdot)\xrightarrow{v} \E N((s,t]\times \cdot)$ in
$M(\Ro)$   and
\begin{equation*}\label{p:poiss}
\lim_{n\to\infty}\Pr(N_n((s,t]\times U)=0)=\Pr(N((s,t]\times U)=0)
\end{equation*}
for every $U\in\mathcal{U}$, then \eqref{eq:1dpc} holds.
\end{enumerate}
\end{lemma}

\section{Limit theorems for stationary sequences}\label{s:sms}

In this section we study limit theorems for stationary sequences of random vectors
in $\st$. Let $\Xa$ be a L\'evy $\alpha$-stable process as in \eqref{eq:le1} with
L\'evy measure $\Pia=\Pi_\alpha$ defined by~\eqref{d:Pia} and let $\Na$ be a
Poisson point process on $[0,\infty)\times \Ro$ with mean measure
$\mathrm{Leb}\times \Pia$.

We assume throughout this section that $\{Z_j\colon j\ge 1\}$ is a strictly
stationary sequence of random vectors in $\st$ such that $Z_1$ is regularly
varying with index $\alpha\in(0,2)$. Then we have
\begin{equation}\label{d:rv1d}
n\Pr(b_n^{-1}Z_1\in\cdot)\xrightarrow{v}\Pia(\cdot),
\end{equation}
where the normalizing constants $b_n$ are such that
\[
\lim_{n\to\infty}n\Pr(|Z_1| > b_n)=1.
\]
We define
\[
X_n(t)=\frac{1}{b_n}\bigl(\sum_{j\le nt}Z_j-tn \E(Z_1I(|Z_1|\le b_n))\bigr)\quad
\text{and}\quad N_n=\sum_{j\ge 1}\epsilon_{(\frac{j}{n},\frac{Z_j}{b_n})}.
\]

\begin{theorem}\label{t:stable} Suppose that $Z_1$ is regularly varying with index
$\alpha\in(0,2)$. Then
 $X_n\xrightarrow{d} \Xa$ in $\D$ if and only if
$N_n\xrightarrow{d} \Na$ in $\ME$ and \eqref{l:rho1} or \eqref{l:rho2} of Theorem~\ref{t:rho} holds.
\end{theorem}
\begin{proof}
Let $X_{n,j}=Z_j/b_n$, $j\ge 1$,  and $ c_n= n\E(X_{n,1} I(|X_{n,1}|\le 1 ))$, $n\ge1$. From~\eqref{d:rv1d} it
follows that, for any $\varepsilon\in(0,1)$,
\begin{equation*}
\lim_{n\to\infty}n\E(X_{n,1} I(\varepsilon < |X_{n,1}|\le 1 )) =\int_{\{x\colon \varepsilon<|x|\le 1\}}x\Pi(dx),
\end{equation*}
which together with $\E(|X_{n,1}| I(|X_{n,1}|\le \varepsilon))\to 0$, as $n\to\infty$, implies that
\[
\lim_{n\to\infty}\sup_{0\le t\le T}\bigl|t\int_{\{x\colon \varepsilon<|x|\le 1\}}x\Pi(dx)-tc_n+[nt]\E(X_{n,1}
I(|X_{n,1}|\le \varepsilon))\bigr|=0,
\]
for all $T>0$ and $\varepsilon\in(0,1)$. Now observe that, by stationarity, condition~\eqref{eq:esmall1d1} holds
for all $\delta>0$ if and only if condition~\eqref{eq:esmall1d} holds for all $T>0$ and $\delta>0$.
Consequently, by Theorem~\ref{th:main} and Remark~\ref{l:iddis}, we obtain $X_n\xrightarrow{d} \Xa$ in $\D$ if
and only if $N_n\xrightarrow{d} \Na$ in $\ME$ and condition~\eqref{eq:esmall1d1} holds for every $\delta>0$.

It remains to show that if $\alpha<1$ then~\eqref{d:rv1d} implies~\eqref{eq:esmall1d1}. By Lemma~\ref{l:max}, we
have
\[
\Pr\bigl(\max_{1\le k\le n}|\sum_{j=1}^{k}(Z_{j}I(|Z_{j}|\le \varepsilon b_n)-\E(Z_{1}I(|Z_{1}|\le \varepsilon
b_n)))|\ge \delta b_n\bigr)\le \frac{2n}{\delta b_n}\E(|\ob|I(|\ob|\le \varepsilon b_n))
\]
for all $n\ge 1$, $\delta,\varepsilon>0$. The rest of the argument is standard. From Karamata's theorem, it
follows that
\[
\E(|\ob|I(|\ob|\le \varepsilon b_n))\sim \frac{\alpha}{1-\alpha}\varepsilon b_n
\Pr(|\ob|\ge \varepsilon b_n)\sim \frac{\alpha}{1-\alpha}(\varepsilon
b_n)^{1-\alpha}L(\varepsilon b_n),
\]
where $L$ is a slowly varying function such that  $n b_n^{-\alpha}L(b_n)\to 1$.
Consequently, we obtain
\[
\limsup_{n\to\infty}\frac{n}{b_n}\E(|\ob|I(|\ob|\le \varepsilon b_n))=\frac{\alpha
}{1-\alpha}\varepsilon^{1-\alpha}
\]
for every $\varepsilon>0$, which completes the proof.
\end{proof}

The next result gives necessary conditions for convergence of point processes to
the Poisson process $\Na$.

\begin{theorem}\label{c:poi} Suppose that $Z_1$ is regularly varying.
If $N_n\xrightarrow{d} \Na$ in $\ME$ then for any $t,\varepsilon>0$ we have
\begin{equation*}
\lim_{n\to\infty}\Pr(\max_{1\le j\le nt}|Z_j|\le \varepsilon
b_n)=e^{-t\Pia(\{x\colon |x|>\varepsilon\})}
\end{equation*}
and
\begin{equation}\label{eq:asneg2}
\lim_{n\to\infty}\Pr(\max_{2\le j\le nt}|Z_j|>  \varepsilon b_n
\bigl||Z_1|>\varepsilon b_n)=1-e^{-t\Pia(\{x\colon |x|>\varepsilon\})}.
\end{equation}
\end{theorem}
The first statement is a consequence of the assumption and the identity
\[
\Pr(N_n((0,t]\times \{x\colon |x|>\varepsilon\})=0)=\Pr(\max_{1\le j\le
nt}|Z_j|\le \varepsilon b_n).
\]
Condition \eqref{eq:asneg2} follows from the next lemma.

\begin{remark}\label{r:asneg}
Observe that the convergence in~\eqref{eq:asneg2} is locally uniform with respect
to $t$. Hence, condition~\eqref{eq:asneg} holds for  every $\varepsilon>0$ and all
sequences $r_n$ such that $r_n=o(n)$.
\end{remark}

\begin{lemma}\label{l:exin}
Let  $\{\xi_j\colon j\ge 1\}$ be a strictly stationary sequence of random
variables. Suppose that $\lambda>0,\theta>0$, and  $u_n$, $n\ge 1$, are such that
\begin{equation}\label{eq:exin}
\lim_{n\to\infty}n\Pr(\xi_1>u_n)=\lambda\quad \text{and}\quad
\lim_{n\to\infty}\Pr(\max_{1\le j\le nt}\xi_j\le u_n)= e^{-\theta\lambda t}
\end{equation}
for all $t>0$. Then
\begin{equation}
\lim_{n\to\infty}\Pr(\max_{2\le j\le nt}\xi_j>u_n\bigl|\xi_1>u_n)= 1-\theta
e^{-\theta\lambda t},\quad t>0.
\end{equation}
\end{lemma}
\begin{proof}
Define for $n\ge 1$, $t\ge 0$
\[F_n(t)=\Pr(\max_{2\le j\le nt}\xi_j>
u_n)
\]
and
\[
G_n(t)=\Pr(\max_{2\le j\le nt}\xi_j> u_n \bigl|\xi_1>u_n),
\]
where we set $F_n(t)=G_n(t)=0$ for $0\le t<2/n$.  Both functions are nondecreasing
and piecewise constant. We first show that
\begin{equation}\label{d:se}
1-G_n(t)=\frac{F_n(t+\frac{1}{n})-F_n(t)}{\Pr(\xi_1>u_n)},\quad t\ge 0, n\ge 1.
\end{equation}
Observe that equality \eqref{d:se} holds for $t\in [0,2/n)$. Let $k\ge 2$ and $t\in [k/n,(k+1)/n)$. We have
\[
1-G_n(t)=\Pr(\max_{2\le j\le k}\xi_j\le u_n \bigl|\xi_1>u_n),
\]
 which leads to
\[
1-G_n(t)=\frac{\Pr(\max_{2\le j\le k}\xi_j\le u_n )-\Pr(\max_{1\le j\le k}\xi_j\le u_n )}{\Pr(\xi_1>u_n)},
\]
and, by stationarity, concludes the proof of \eqref{d:se}.

To complete the proof it suffices to show that
\[
n(F_n(t+\frac{1}{n})-F_n(t))\to \theta\lambda e^{-\theta\lambda t}.
\]
We proceed similarly to \cite[pp. 2047-2048]{haydnlacroix05}. Define piecewise
linear functions $\widetilde{F}_n$ by $\widetilde{F}_n(t)=F_n(t)$ for $t=k/n$ and
$\widetilde{F}_n$ linear on $[k/n,(k+1)/n]$, $k\ge 0$. Then the right-hand
derivative $\widetilde{F}_n'(t+)$ at every point $t$ is given by
\[
\widetilde{F}_n'(t+)=n(F_n(t+\frac{1}{n})-F_n(t)).
\]
Note that $F_n(t)\to 1-e^{-\theta\lambda t}$, as $n\to\infty$, for all $t\ge 0$,
and, by \eqref{d:se}, we have
\[
\sup_{t\ge 0}|\widetilde{F}_n(t)-F_n(t)|\le \Pr(\xi_1>u_n)\to 0.
\]
Since the functions $\widetilde{F}_n$ are concave, we obtain
$\widetilde{F}_n'(t+)\to \theta\lambda e^{-\theta\lambda t}$ for all $t$.
\end{proof}

\begin{remark}\label{r:exi}
 Note that the constant $\theta$ in~\eqref{eq:exin} might be referred to
as the extremal index of the sequence $\{\xi_j\colon j\ge 1\}$; see e.g.~\cite{leadbetterrootzen88,obrien87}
and~\cite[Chapter 8.1]{embrechts97}  for the definition and properties. In particular, if the $\xi_j$ are
i.i.d.~then $\theta=1$.  Dependent random variables have the extremal index equal to $1$ when they satisfy the
extreme mixing conditions $D(u_n)$ and $D'(u_n)$.  This will be also the case for the sequence $\{|Z_j|\colon
j\ge 1\}$ in Theorem~\ref{t:stable}.
\end{remark}

We now provide sufficient conditions for
 convergence to Poisson processes for strongly mixing sequences.
The mixing condition $D^*$ of Davis and Resnick~\cite[p.~47]{davisresnick88} is implied by strong mixing. Hence,
from \cite[Theorem 2.1]{davisresnick88} it follows that if the local dependence condition~\eqref{eq:D'} of
Davis~\cite{davis83}, then $N_n\xrightarrow{d} \Na$ in $\ME$. Although condition~\eqref{eq:D'} is sufficient for
Poisson convergence, it is not necessary. We now prove that our condition \textbf{LD}$(\phi_0)$ from
Theorem~\ref{t:rho}, which is necessary, is also sufficient.

\begin{theorem}\label{th:pm}
Suppose that the sequence $\{Z_j\colon j\ge 1\}$  is strongly mixing and $Z_1$ is regularly varying. If
condition~\textbf{LD}$(\phi_0)$ holds then $N_n\xrightarrow{d} \Na$ in $\ME$.
\end{theorem}
\begin{proof}
By part~\eqref{l:pcl1} of Lemma~\ref{l:pc}, we have $N_n\xrightarrow{d} \Na$ if and only if for any $t>s\ge 0$
condition~\eqref{eq:1dpc} holds, since $\{Z_j\colon j\ge 1\}$ is strongly mixing. For any $f\in C_K^{+}(\Ro)$
and $t>s\ge 0$ we have, by stationarity,
\[
\bigl|\E(e^{-\sum_{ns<j\le nt}f(b_n^{-1}Z_j)})-\E(e^{-\sum_{0<j\le n(t-s)}f(b_n^{-1}Z_j)})\bigr|\le
\E(1-e^{-f(b_n^{-1}Z_1)}),
\]
which converges to $0$ as $n\to\infty$, since $f(b_n^{-1}Z_1)\to 0$ a.s. Consequently, we have
$N_n\xrightarrow{d} \Na$ if and only if for any $t>0$
\begin{equation*}
N_n((0,t]\times \cdot)\xrightarrow{d} \Na((0,t]\times \cdot)\quad \text{in } \MEE.
\end{equation*}
Let $t>0$. From~~\eqref{d:rv1d} it follows that $\E(N_n((0,t]\times \cdot))\xrightarrow{v}\E(N((0,t]\times
\cdot))$. Hence, it suffices to show, by~part~\eqref{l:pcl2} of Lemma~\ref{l:pc}, that
\begin{equation}\label{eq:estpoi}
\lim_{n\to\infty}\Pr(N_n((0,t]\times U)=0)=e^{-t\Pia(U)}
\end{equation}
for every finite union $U$ of rectangles with $\Pia(U)<\infty$ and $\Pia(\partial
U)=0$.

Let $\varepsilon>0$ be such that $U\subset\{x\colon |x|>\varepsilon\}$. Take
$r_n,l_n$ as in~\eqref{eq:asneg0} and~\eqref{eq:asneg}. Since $l_n=o(r_n)$ we may
assume that $l_n<r_n$. Let the integers $k_n,s_n$ be given by the Euclidean
division of $[nt]$ by $r_n$, $[nt]=k_nr_n+s_n$ and $0\le s_n<r_n$. We will prove
the following two statements
\begin{equation}\label{eq:pom1}
|\Pr(N_n((0,t]\times U)=0)-\bigl(1-\Pr(\bigcup_{j=1}^{r_n-l_n} \{Z_j\in b_n
U\})\bigr)^{k_n}|\to 0
\end{equation}
and
\begin{equation}\label{eq:pom2}
\frac{\Pr(\bigcup_{j=1}^{r_n-l_n} \{Z_j\in b_n U\})}{r_n\Pr(Z_1\in b_n U)}\to 1.
\end{equation}
Since $k_nr_n\Pr(A_1)\to t\Pia(U)$, conditions~\eqref{eq:pom1} and~\eqref{eq:pom2}
imply~\eqref{eq:estpoi}.

To prove~\eqref{eq:pom1} we use the standard big-little block technique.  Write
$A_j=\{Z_j\in b_n U\}$, $j\ge 1$, and observe that
\[
|\Pr(\bigcap_{j=1}^{[nt]}A_j^c)-\Pr(\bigcap_{j=1}^{k_nr_n}A_j^c)|\le r_n\Pr(A_1).
\] Let us divide
the integers $1,\ldots,k_nr_n$, into blocks of the form
\[
I_j=\{(j-1)r_n+1,\ldots,jr_n-l_n\},\quad
I_j^*=\{(j-1)r_n-l_n+1,\ldots,jr_n\},\quad j=1,\ldots,k_n.
\]
We have
\[
|\Pr(\bigcap_{j=1}^{k_nr_n}A_j^c)-\Pr(\bigcap_{j=1}^{k_n}\bigcap_{i\in
I_j}A_i^c)|\le \Pr(\bigcup_{j=1}^{k_n}\bigcup_{i\in I_j^{*}}A_i^c)\le
k_nl_n\Pr(A_1)
\]
and, by strong mixing,
\[
|\Pr(\bigcap_{j=1}^{k_n}\bigcap_{i\in I_j}A_i^c)-\Pr(\bigcap_{i\in
I_1}A_i^c)^{k_n}|\le (k_n-1)\phi_0(l_n).
\]
Summarizing
\begin{equation*}\label{eq:pom}
|\Pr(\bigcap_{j=1}^{[nt]}A_j^c)-\Pr(\bigcap_{j\in I_1}A_j^c)^{k_n}|\le (r_n +k_nl_n)\Pr(A_1)+(k_n-1)\phi_0(l_n)
\end{equation*}
and, by the choice of the sequences, the right-hand side in the last inequality
goes to $0$ as $n\to\infty$, which completes the proof of~\eqref{eq:pom1}. Now
observe that
\[
|\Pr(\bigcup_{j=1}^{r_n}A_j)-\Pr(\bigcup_{j=1}^{r_n-l_n}A_j)|\le 2l_n\Pr(A_1).
\]
Hence, it remains to show that
\[
\frac{\Pr(\bigcup_{j=1}^{r_n}A_j)}{r_n\Pr(A_1)}\to 1.
\]
Since $A_j=\{Z_j\in b_n U\} \subset \{|Z_j|>\varepsilon b_n\}$, we have
\[
\Pr(\bigcup_{j=2}^{r_n}A_j\bigl|A_1)\le \Pr(\max_{2\le j\le r_n}|Z_j|> \varepsilon
b_n\bigl||Z_1|>\varepsilon b_n)\frac{\Pr(|Z_1|> \varepsilon b_n)}{\Pr(Z_1\in
b_nU)},
\]
which shows that the left-hand side in the last inequality goes to $0$ as
$n\to\infty$. Consequently,
\begin{equation}\label{eq:pom3}
\Pr(\bigcup_{j=3}^{r_n}A_j\bigl|A_1)\to 0\quad \text{and}\quad
\Pr(A_2\bigl|A_1)\to 0.
\end{equation}
We have
\[
\begin{split}
\Pr(\bigcup_{j=1}^{r_n}A_j)&=\sum_{j=1}^{r_n-2}\Pr(A_j\cap A_{j+1}^c\cap
\bigcap_{i=j+2}^{r_n}A_i^c)+\Pr(A_{r_n-1}\cap A_{r_n}^c)+\Pr(A_{r_n})
\end{split}
\]
and thus, by stationarity,
\[
\begin{split}
|r_n\Pr(A_1\cap A_2^c)-\Pr(\bigcup_{j=1}^{r_n}A_j)|&\le \sum_{j=1}^{r_n-2}\Pr(A_j\cap A_{j+1}^c\cap
\bigcup_{i=j+2}^{r_n}A_i)+\Pr(A_1\cap
A_2)\\
&\le r_n\Pr(A_1\cap \bigcup_{i=3}^{r_n}A_i)+\Pr(A_1\cap A_2),
\end{split}
\]
since $\Pr(A_1\cap A_2^c)=\Pr(A_j\cap A_{j+1}^c)$ and $\Pr(A_j\cap A_{j+1}^c\cap
\bigcup_{i=j+2}^{r_n}A_i)=\Pr(A_1\cap A_2^c \cap \bigcup_{i=3}^{r_n +1-j}A_i)$ for each $j=1,\ldots,r_n-2$,
which completes the proof by \eqref{eq:pom3}.
\end{proof}

Theorem~\ref{t:rho} is  a direct consequence of Theorems~\ref{t:stable},~\ref{c:poi}, and~\ref{th:pm}. For the
proof of Corollary~\ref{c:kobus} we need the following result of Novak~\cite[Corollary 2.2]{novak96}.

\begin{lemma}\label{l:novak} Let  $\{\xi_j\colon j\ge 1\}$ be a strictly stationary and uniformly mixing sequence of random
variables. If the sequence $u_n$  is such that
\[
0<\liminf_{n\to\infty} n\Pr(\xi_1>u_n)\le
\limsup_{n\to\infty}n\Pr(\xi_1>u_n)<\infty
\]
and, for every $j\ge 2$,
\[
\lim_{n\to\infty}\Pr(\xi_j>u_n\bigl|\xi_1>u_n)=0,
\]
 then
\[
\Pr(\max_{1\le j\le n}\xi_j\le u_n)-\exp(-n\Pr(\xi_1>u_n))\to 0.
\]
\end{lemma}

\begin{proof}[Proof of Corollary~\ref{c:kobus}] Let $\varepsilon>0$ and $t>0$. Define $u_n=\varepsilon b_{[n/t]}$,
$n\ge 1$.  Set $\lambda=\Pia(\{x\colon |x|>\varepsilon\})$  and observe that
$n\Pr(|Z_1|>\varepsilon b_n)\to \lambda$ and $n\Pr(|Z_1|>u_n)\to \lambda t$, as
$n\to\infty$. From \eqref{d:watson} it follows that the sequences $u_n$ and
$\{|Z_j|\colon j\ge 1\}$ satisfy all assumptions of Lemma~\ref{l:novak}. Hence,
\[
\lim_{n\to\infty}\Pr(\max_{1\le j\le n}|Z_j|\le u_n)= e^{-\lambda t},
\]
and, consequently,
\[
\lim_{n\to\infty}\Pr(\max_{1\le j\le nt}|Z_j|\le \varepsilon b_n)=e^{-\lambda t}.
\]
From Lemma~\ref{l:exin} it follows that
\[
\lim_{n\to\infty}\Pr(\max_{2\le j\le nt}|Z_j|>\varepsilon b_n\bigl||Z_1|>\varepsilon b_n)=1-e^{-\lambda t},
\]
which implies condition~\eqref{eq:asneg}, by Remark~\ref{r:asneg}, and concludes the proof.
\end{proof}

\begin{lemma}\label{lemat}
Suppose that the maximal correlation coefficient $\rho(n)=\rho(\F_{1}^{1},\F_{n+1}^{\infty})$ of the sequence
$\{Z_j\colon j\ge 1\}$ satisfies $\sum_{j}\rho(2^j)<\infty$. If $Z_1$ is regularly varying with index
$\alpha\in[1,2)$ then condition~\eqref{eq:esmall1d1} holds for all $\delta>0$.
\end{lemma}
\begin{proof}
First observe that it suffices to show that for every $\delta>0$ there exists a constant $C>0$ such that for any
$\varepsilon>0$ and $n\ge 1$
\[
\Pr\bigl(\max_{1\le k\le n}|\sum_{j=1}^{k}(Z_{j}I(|Z_{j}|\le \varepsilon b_n)-\E(Z_{1}I(|Z_{1}|\le \varepsilon
b_n))|\ge \delta b_n\bigr)\le  \frac{C n}{b_n^{2}}\E(|Z_{1}|^2I(|Z_{1}|\le \varepsilon b_n)),
\]
since $|Z_1|$ is regularly varying with index $\alpha<2$ and  $n b_n^{-2}\E(|Z_{1}|^2I(|Z_{1}|\le \varepsilon
b_n))\to \frac{\alpha}{2-\alpha}\varepsilon ^{2-\alpha}$, by Karamata's theorem. When $d=1$ then this type of
bound follows from the $L^2$-maximal inequality from~\cite{shao95}. We now outline how to get a similar bound in
the multivariate case.

Write $Z_{j}(a)=Z_{j}I(|Z_{j}|\le a)-\E(Z_{j}I(|Z_{j}|\le a))$ for $a>0$, $j\ge 1$. For every $a$ and $j$  the
random vector $Z_{j}(a)$ has zero mean and is bounded.
 The proof of
the $L^2$-maximal inequality for stationary sequences of random variables as given
in~\cite[pp.~544-555]{peligradutevwu} still works for random vectors and we can deduce the following
\[
\E(\max_{1\le k\le n}|S_k|^2)\le 2 n \Bigl(2\|Z_1(a)\|_2+ 4 \sum_{j=0}^{[\log_2n]}
2^{-j/2}\|\E(S_{2^{j}}|Z_1(a))\|_2\Bigr)^2,
\]
where $S_k=\sum_{j=2}^{k+1} Z_j(a)$ for $k\ge 1$ and $\|Y\|_2=\sqrt{\E\langle Y,Y\rangle}$. By Chebyshev's
inequality,  it remains to show that there exists a constant $C_1$ such that for any $a>0$
\[
\sum_{j=0}^\infty 2^{-j/2}\|\E(S_{2^j}|Z_1(a))\|_2\le C_1 \|Z_1(a)\|_2,
\]
since $\|Z_1(a)\|_2^2\le 2 \E (|Z_1|^2I(|Z_1|\le a))$. We have
\[
\|\E(S_{2n}|Z_1(a))\|_2\le \|\E(S_{n}|Z_1(a))\|_2+ \|\E(S_{2n}-S_n|Z_1(a))\|_2,\quad n\ge 1
\]
and, by using \cite[Theorem 4.2]{bradleybryc85},
\[
\begin{split}
 \|\E(S_{2n}-S_n|Z_1(a))\|_2^2&=\E\langle S_{2n}-S_n,\E(S_{2n}-S_n|Z_1(a))\rangle\\
 &\le \rho(n)\|S_{2n}-S_n\|_2 \|\E(S_{2n}-S_n|Z_1(a))\|_2,
 \end{split}
\]
which implies
\[
 \|\E(S_{2n}-S_n|Z_1(a))\|_2\le\rho(n)\|S_{2n}-S_n\|_2=\rho(n)\|S_n\|_2.
\]
The proof of \cite[Lemma 3.4]{peligrad82} extends directly to random vectors. Thus, there exits a constant $C_2$
such that for every $a>0$ and $n\ge 1$ we have
\[
\|S_n\|_2\le C_2 \sqrt{n}\|Z_1(a)\|_2,
\]
which gives, as in \cite[Lemma 1]{peligradutevwu}, the following estimate
\[
\sum_{j=0}^\infty 2^{-j/2}\|\E(S_{2^j}|Z_1(a))\|_2\le 4C_2 \|Z_1(a)\|_2 \sum_{j=0}^\infty \rho(2^j),
\]
and completes the proof.
\end{proof}

\begin{remark}
Note that since the sequence $\rho(n)$, $n\ge 1$, is nonincreasing we have
\[
\sum_{j=0}^\infty \rho(2^j)<\infty\quad \text{if and only if}\quad \sum_{n=1}^\infty\frac{\rho(n)}{n}<\infty.
\]
\end{remark}

Corollary~\ref{c:mdep} follows from Theorems~\ref{t:rho},~\ref{th:main}, Lemma~\ref{lemat}, and \cite[Theorem
1]{hudson89} or \cite[Theorem 1.1]{kobus95}.

\section{Proofs of Theorems~\ref{th:main} and~\ref{p:marg}}\label{s:proof}

\begin{proof}[Proof of Theorem~\ref{th:main}] Since the L\'evy process has no fixed points of discontinuity, it
follows that $X_n\xrightarrow{d} X$ in $\D$  if and only if $X_n\xrightarrow{d} X$
 in $\DT$ for any $T>0$.

First assume that \eqref{eq:pce} and \eqref{ce:assmall} hold. For the proof of \eqref{d:con} we adapt the
arguments of \cite[Section 4]{durrettresnick78} (see also \cite[Section 7.2]{resnick07}). Let us define
\[
\eta_\varepsilon^{(1)}(t)=\int_{[0,t]\times\{x:|x|> \varepsilon\}}x N(ds, dx)
\]
for $\varepsilon >0$ and for every $\varepsilon\in[0,1)$
\[
\eta_{\varepsilon}^{(2)}(t)=\int_{[0,t]\times\{x:\varepsilon<|x|\le
1\}}x\bigl(N(ds, dx)-ds\Pi(dx)\bigr),\quad t\ge 0.
\]
By the L\'evy-It\^{o} integral representation, we can rewrite $X$ almost surely as
\begin{equation}\label{e:levy}
 \eta(t)=\eta_1^{(1)}(t)+\lim_{\varepsilon\downarrow
 0}\eta_\varepsilon^{(2)}(t).
\end{equation}
The terms in \eqref{e:levy} are independent and the convergence in the last term
is a.s. and uniform in $t$ on any bounded interval. Hence, we obtain
\begin{equation}\label{e:levy0}
\eta_1^{(1)}+\eta_\varepsilon^{(2)}\xrightarrow{d}\eta\quad\text{in }\D,\quad
\text{as }\varepsilon \to 0.
\end{equation}
From~\eqref{eq:pce} and the continuous mapping theorem it follows that
\[
\eta_{n,\varepsilon}^{(1)}\xrightarrow{d} \eta_{\varepsilon}^{(1)}\quad\text{in
}\D, \quad \text{as }n \to \infty,
\]
for all $\varepsilon\in(0,1)$ such that $\Pi(\{x:|x|=\varepsilon\})=0$, where
\[
 \eta_{n,\varepsilon}^{(1)}(t):=\sum_{j\le
nt} \eta_{n,j}I(|\eta_{n,j}|>\varepsilon),\quad t\ge 0,
\]
since the mapping $R_{0,\varepsilon}\colon \ME\to\D$ defined by
\[
R_{0,\varepsilon}(m)(t)=\int_{[0,t]\times\{x:|x|> \varepsilon\}}x\, m(ds, dx), \quad m\in\ME,
\]
is a.s.~continuous with respect to the distribution of the Poisson point process $N$ for all such
$\varepsilon\in(0,1)$ (see e.g. \cite[p.84]{resnick86} or \cite[Section 7.2]{resnick07}).
 Hence, for
\[
\eta_{n,\varepsilon}(t)=\eta_{n,\varepsilon}^{(1)}(t)-t\int_{\{x:\varepsilon
<|x|\le 1\} }x\Pi(dx),\quad t\ge 0,
\]
we obtain $\eta_{n,\varepsilon}\xrightarrow{d}\eta_1^{(1)}+\eta_\varepsilon^{(2)}$ in $\D$. The function
$\varepsilon\mapsto \Pi(\{x:|x|>\varepsilon\})$ is monotonic. Therefore, we may chose a sequence
$\varepsilon_k\in(0,1)$ such that $\Pi(\{x:|x|=\varepsilon_k\})=0$ and $\varepsilon_k\downarrow 0$ From
\eqref{e:levy0} and the converging together theorem \cite[Theorem 4.2]{billingsley68}, it suffices to show that,
for any $\delta>0$,
\[
\lim_{k\to\infty}\limsup_{n\to\infty}\Pr\bigl(d_{T}(X_n,X_{n,\varepsilon_k})\ge \delta\bigr)=0.
\]
This is a consequence of \eqref{ce:assmall}, since
\[
\eta_n(t)- \eta_{n,\varepsilon}(t)=\sum_{j\le nt}\eta_{n,j}I(|\eta_{n,j}|\le
\varepsilon)-tc_n+t\int_{\{x:\varepsilon <|x|\le 1\} }x\Pi(dx)
\]
and the Skorohod metric $d_T$ on $\DT$ is bounded above by the uniform metric on $\DT$.

Now assume that \eqref{d:con} holds. To prove \eqref{eq:pce} it suffices to show
that for every $f\in C_K^+([0,\infty)\times \Ro)$ we have
\[
\E(e^{-N_n(f)})\to \E(e^{-N(f)}).
\]
Let $U(X)=\{r>0:\Pr(|\Delta \eta(t)|=r\text{ for some }t>0)>0\}$. The set
$U(\eta)$ is at most countable \cite[Lemma VI.3.12]{jacodshiryaev03}. Let
$\varepsilon\not\in U(\eta)$ and $E_{T,\varepsilon}=[0,T]\times \{x:
|x|>\varepsilon\}$ for $T>0$. Define the mapping $R_{1,\varepsilon}\colon \D\to
M_p(E_{T,\varepsilon})$ by
\[
R_{1,\varepsilon} \psi=\sum_{\{t\le T:
|\Delta\psi(t)|>\varepsilon\}}\epsilon_{(t,\Delta\psi(t))}.
\]
Since for every $\psi\in\D$ the set $\{t\le T: |\Delta\psi(t)|>\varepsilon\}$ is
finite, the mapping $R_{1,\varepsilon}$ is well defined. Moreover,
$R_{1,\varepsilon}$ is continuous at all $\psi$ such that $\varepsilon\not\in
\{r>0: |\Delta\psi(t)|=r\text{ for some } t>0\}$ (see e.g.~\cite[Section
6.2]{jacodshiryaev03}). Hence, the mapping $R_{1,\varepsilon}$ is almost surely
continuous with respect to the distribution of $X$. From the continuous mapping
theorem it follows that
\[
R_{1,\varepsilon}X_n\xrightarrow{d}R_{1,\varepsilon}X\quad \text{in
}M_p(E_{T,\varepsilon}).
\]
Thus
\[
\E(e^{-R_{1,\varepsilon}X_n(f)})\to \E(e^{-R_{1,\varepsilon}X (f)})\quad \text{for
all } f\in C_K^+(E_{T,\varepsilon}).
\]
Observe that we have $\Delta X_n(t)\neq 0$ if and only if $t=j/n$ for some $j$.
Since for every $f\in C_K^+([0,\infty)\times \Ro)$ we can find $T>0$ and
$\varepsilon>0$ such that the support of $f$ is contained in $E_{T,\varepsilon}$,
we obtain
\[
\E(e^{-N_{n}(f)})=\E(e^{-R_{1,\varepsilon}X_n(f)})\to \E(e^{-R_{1,\varepsilon}X
(f)})=\E(e^{-N(f)}),
\]
by the definition of $N_n$ in \eqref{eq:pce} and that of $N$ in~\eqref{d:pp}, which completes the proof
of~\eqref{eq:pce}.

To prove \eqref{ce:assmall}, we first show that
\[
X_n-\sum_{s\le \cdot}\Delta X_n(s)I(|\Delta
X_n(s)|>\varepsilon)\xrightarrow{d}X-\sum_{s\le \cdot}\Delta X(s)I(|\Delta
X(s)|>\varepsilon)\quad\text{in }\D
\]
for all $\varepsilon>0$ such that $\Pi(\{x:|x|=\varepsilon\})=0$. Define the
mapping $R_{2,\varepsilon}\colon \D\to \D$ by
\[
R_{2,\varepsilon}\psi(t)=\psi(t)-\sum_{s\le
t}\Delta\psi(s)I(|\Delta\psi(s)|>\varepsilon),\quad t\ge 0.
\]  By \cite[Proposition VI.2.7]{jacodshiryaev03}, $R_{2,\varepsilon}$ is
continuous at $\psi$ if $\varepsilon\not\in\{r\colon |\Delta\psi(t)|=r \text{ for
some }t>0\}$. Observe that $\{\varepsilon>0\colon \Pi(\{x\colon
|x|=\varepsilon\})=0\}\subseteq \mathbb{R}_+\setminus U(X)$, thus the claim
follows from the continuous mapping theorem. For $\varepsilon <1$ and $\psi\in\D$
define
\[
R_\varepsilon \psi(t):=R_{2,\varepsilon}\psi(t)+t\int_{\{x\colon \varepsilon
<|x|\le 1\} }x\Pi(dx).
\]
We have
\[
R_\varepsilon X_n(t)=\sum_{j\le nt}X_{n,j}I(|X_{n,j}|\le
\varepsilon)-t(c_n-\int_{\{x\colon \varepsilon <|x|\le 1\} }x\Pi(dx))
\]
and
\[
R_\varepsilon X(t)=X(t)-X_{1}^{(1)}(t)-X_\varepsilon^{(2)}(t),\quad t\ge 0.
\]
The set $F_\delta=\{\psi\in\D\colon \sup_{0\le t\le T}|\psi(t)|\ge \delta\}$ is
closed in $\D$. Since $R_\varepsilon X_n\xrightarrow{d}R_\varepsilon X$, we obtain
\[
\limsup_{n\to\infty}\Pr(R_\varepsilon X_n\in F_\delta)\le \Pr(R_\varepsilon X\in
F_\delta),
\]
 by Portmanteau's
theorem. From \eqref{e:levy} it follows that
\[
\lim_{\varepsilon\to 0}\Pr(R_\varepsilon X\in F_\delta)=0
\]
 which completes the proof of \eqref{ce:assmall}.
\end{proof}

\begin{proof}[Proof of Theorem~\ref{p:marg}]
Let $\varepsilon\in(0,1)$ be such that $\Pi(\{x:|x|=\varepsilon\})=0$ and let
$f_\varepsilon(x)=xI(|x|>\varepsilon)$. Since $\Pi(1\wedge
|f_\varepsilon|)<\infty$, the random vector $N'(f_\varepsilon)$ has the
characteristic function of the form (see e.g.~\cite[Lemma 12.2]{kallenberg02})
\[
\E(e^{i\langle u,N'(f_\varepsilon)\rangle})=\exp(\int (e^{i\langle u,
x\rangle}-1)I(|x|>\varepsilon)\Pi(dx)).
\]
From \eqref{eq:pce1d} and the continuous mapping theorem it follows that
\begin{equation*}\label{e:creps}
N_n'(f_\varepsilon)\xrightarrow{d} N'(f_\varepsilon)\quad\text{in }\st.
\end{equation*}
With the notation as in the proof of Theorem~\ref{th:main}, observe that
$\eta_1^{(1)}+\eta_\varepsilon^{(2)}$  is a L\'evy process such that
\[
\E  (e^{i\langle u,\eta_1^{(1)}(1)+\eta_\varepsilon^{(2)}(1) \rangle})=\exp
\bigl(\int \bigl(e^{i\langle u, x\rangle}-1-i\langle u, x\rangle
I(\varepsilon<|x|\le 1)\bigr)\Pi(dx)\bigr)
\]
and, by \eqref{e:levy0},
\[
\eta_1^{(1)}(1)+\eta_\varepsilon^{(2)}(1)\xrightarrow{d}\eta(1)\quad \text{in }
\st.
\]
Since $N'(f_\varepsilon)-\int_{\{x:\varepsilon <|x|\le 1\} }x\Pi(dx)$ has the same
distribution as $\eta_1^{(1)}(1)+\eta_\varepsilon^{(2)}(1)$, the result follows
from \eqref{ce:assmall1d} and \cite[Theorem 4.2]{billingsley68}.
\end{proof}

\section*{Acknowledgments} This work was supported by
Polish MNiSW grant N N201 0211 33 and by EPSCR grant  EP/F031807/1 \emph{Anomalous Diffusion in Deterministic
Systems}. Parts of this work were completed while the author was visiting the University of Surrey, Guildford,
United Kingdom, and the hospitality of the Department of Mathematics is gratefully acknowledged. We thank the
referees for constructive remarks which led to an improved presentation of the paper.

\end{document}